\documentclass{article}
\title{About calculation of the Hankel transform
using preliminary wavelet transform}
\author{E.B. Postnikov}
\date{Theoretical Physics Department,
Kursk State Pedagogical University, \\ Radischsheva st. 33, Kursk,
Russia, 305000 \\ e-mail: postnicov@mail.ru \\   }
\usepackage{a4,graphicx}
\begin{document}
\maketitle
 MSC-class: 44A15, 65R10
\begin{abstract}
The purpose of this paper is to present an algorithm for
evaluating Hankel transform of the null and the first kind. The
result is the exact analytical representation as the series of the
Bessel and Struve functions multiplied to the wavelet coefficients
of the input function. Numerical evaluation of the test function
with known analytical Hankel transform illustrates proposed
algorithm.
\end{abstract}

       The Hankel transform is a very useful instrument in a
 wide range of physical problems which have an axial symmetry~\cite{MMF} .
The influence of the Laplasian on a function in a cylindrical
coordinates is equal to the product
 of the parameter of the transformation squared and the
  transform of the function:

\begin{equation}
\label{eq1}
\begin{array}{l}
 \left( {\frac{{d^2 }}{{dr^2 }} + \frac{1}{r}\frac{d}{{dr}}} \right)f(r) \leftrightarrow  - p^2 F_0 (p) \\
 \left( {\frac{{d^2 }}{{dr^2 }} + \frac{1}{r}\frac{d}{{dr}} - \frac{1}{{r^2 }}} \right)f(r) \leftrightarrow  - p^2 F_1 (p) \\
 \end{array}
\end{equation}

    The Hankel transform of the null $(n = 0)$
and the first $(n = 1)$ kind are represented as
\begin{equation}
\label{eq2}
\begin{array}{l}
 F_n (p) = \int\limits_0^\infty  {f(r)J_n (pr)rdr,}  \\
 f_n (p) = \int\limits_0^\infty  {F(p)J_n (pr)pdp.}  \\
 \end{array}
\end{equation}

Besides that integrals like (2) are connected with the problems of
 geophysics and cosmology,  for example~\cite{Geo, Astro}.

       But practical calculation of direct and inverse
Hankel transform is connected with two problems.
 The first problem is based on the fact that not every transform
in the real physical situation has analytical expression for result
 of inverse Hankel transform. The second one is determination of
functions as a set of their values for numerical calculations.
Large bibliography on those issues can be found in~\cite{LO} The
classical trapezoidal rule, Cotes rule and other rules connected
with replacement of integrand by sequence of polynoms have high
accuracy if integrand is a smooth function. But $ f(r)J_n (pr)r $
(or $F_p (p)J_n (pr)p$)
              is a quick oscillating function if $r$ (or $p$) is large.
 There are two general methods of the
effective calculation in this area. The first is the Fast Hankel
transform ~\cite{FHT}. The specification of that method is
transforming the function to the logarithmical space and fast
Fourier transform in that space. This method needs a smoothing of
the function in log-space. The second method is based on the
separation of the integrand into product of slowly varying
component and a rapidly oscillating Bessel function~\cite{Filon}.
 But it needs the smoothness
of the slow component for its approximation by low-order polynoms.

    The goal of this article is to apply wavelet transform with
 Haar bases to (\ref{eq2}).

     The both direct and inverse transforms (\ref{eq2})
  are  symmetric. Let us consider only one of them,
for example, direct transform. Let's denote $f(r)r$ as $g(r)$.
 Then Hankel transform is
\begin{equation}
\label{eq3}
 F_{0,1} (p) = \int\limits_0^\infty {g(r)J_{0,1}
(pr)dr}
\end{equation}

    The expansion $g(r) \in L^2 (R)$
into wavelet series with the Haar bases is~\cite{Chui} :

\begin{equation}
\label{eq4}
g(r) = \sum\limits_{k = 0}^\infty  {c_{0k} \varphi _k (r)
 + \sum\limits_{j = 0}^\infty  {\sum\limits_{k = 0}^\infty  {d_{jk}
\psi _{jk}^{} (r)} } },
\end{equation}

       \[
       \varphi _{0k} (r) = \varphi ^H (r - k),
                 \psi _{jk} = 2^{j/2} \psi ^H (2^j r - k)
            ,
\]
 \[
    \varphi ^H (t) = \left\{ \begin{array}{l}
        1,t \in (0,1) \\
        0,t \notin (0,1) \\
        \end{array} \right. ,
              \psi ^H (t) = \left\{ {\begin{array}{*{20}c}
          {1,}  \\
          { - 1,}  \\
          {0,}  \\
       \end{array}} \right.\begin{array}{*{20}c}
          {t \in (0,1/2)}  \\
          {t \in (1/2,0)}  \\
          {t \notin (0,1)}  \\
       \end{array}
           .\]

    After substitution  (\ref{eq4}) into  (\ref{eq3}) one has
       \[
       F_{0,1} (p) = \sum\limits_{k = 0}^\infty  {c_{0k}
\int\limits_0^\infty  {\varphi _k (r)J{}_{0,1}(pr)dr}
  + \sum\limits_{j = 0}^\infty  {\sum\limits_{k = 0}^\infty
{d_{jk} \int\limits_0^\infty  {\psi _{jk}^{} (x)J_{0,1} (pr)dr} } } }
.       \]

    Making use integrals from~\cite{SF} we have as a result
    \begin{equation}
\label{eq5}
\begin{array}{l}
F_0(p) =
 \sum\limits_{k \in Z} c_{0k} 
 \left[
(k + 1)J_0 \left( p(k + 1) \right) - kJ_0 \left( pk \right)
\right. +
\\ \\
\frac{\pi}{2}
\left[
(k+1)D\left( p(k + 1) \right)-kD\left( pk \right)
\right]
  + \\ \\
\sum\limits_{j = 0}^\infty  \sum\limits_{k \in Z} d_{jk}
\left[
2(k + \textstyle{1 \over 2} )
J_0 \left( p(k + \textstyle{1 \over 2}) \right)-
(k + 1) J_0 \left( p(k + 1) \right)
- kJ_0 \left( pk \right)-
\right.\\ \\
\frac{\pi}{2}
\left[
2(k + \textstyle{1 \over 2} )
D\left( p(k +\textstyle{1 \over 2} ) \right)-
(k+1)D\left( p(k + 1) \right)-kD\left( pk \right)
\right]
\end{array}
\end{equation}

       \begin{equation}
\label{eq6}
       \begin{array}{l}
        F_1 (p) = \frac{1}{p}\left\{ {\sum\limits_{k \in Z} {c_{0k} \left[ {J_0 \left( {pk} \right) - J_0 \left( {p(k + 1)} \right)} \right] + } } \right. \\
        \left. {\sum\limits_{j = 0}^\infty  {\sum\limits_{k \in Z} {d_{jk} \left[ {2J_0 \left( {p(k + {\textstyle{1 \over 2}})2^{ - j} } \right) - J_0 \left( {p(k + 1)2^{ - j} } \right) - J_0 \left( {pk2^{ - j} } \right)} \right]} } } \right\} \\
        \end{array},
       \end{equation}

       where $D(\xi)=H_0(\xi)J_1(\xi)-H_1(\xi)J_0(\xi)$ and
 $H_{0,1}$
       is Struve function of the null and the first kind.

The most sufficient result is that equations  (\ref{eq5}) and  (\ref{eq6}) are
 exact. They can be used in any analytical expressions.
Especially it's useful for Hankel transform of the first kind
because
 (\ref{eq6})
 contains only a combination of Bessel functions, and one can use such
their
 properties as orthogonality etc.
       The coefficients
       $  c_{0k} $ means average value of $ g(r) $ at the range
       $
       [k,k + 1]
       $
       :
       \[
       c_{ok}  = \int\limits_k^{k + 1} {g(r)dr}
       \]

    The detail coefficients are
       \[
       d_{jk}  = 2^{j/2} \left\{ {\int\limits_{2^{ - j} k}^{2^{ - j} \left( {k + \frac{1}{2}} \right)} {g(r)dr}  - \int\limits_{2^{ - j} \left( {k + \frac{1}{2}} \right)}^{2^{ - j} \left( {k + 1} \right)} {g(r)dr} } \right\}
      . \]

The formulas (\ref{eq5}) and (\ref{eq6}) allow us to get full
analytical solution if integrals above have close form solution.
      In the opposite case the solution must be numerical but this method
provides an effective algorithm for that.
        It's obvious that
       $
       d_{jk}
       $
       decrease very quickly if
       $
       g(r)
       $
       is a smooth function. One can practically use
       $
       d_{jk}  > \varepsilon
       $, where
       $
       \varepsilon
       $
       is small. But if
       $
       g(r)
       $
has steps, sharp vertices or discontinues then the detail coefficients
concentrate around these points and one can appropriate they are
 equal to the zero in other areas.

Let us consider for example a function with known analytical
Hankel transform
\begin{equation}
\label{eq7}
\int\limits_0^\infty{e^{-a^2 r^2}rJ_1 (pr)rdr=
\frac{p}{4a^4}e^{-p^2/4a^2}.}
\end{equation}
The approximation and detail coefficients  may be calculated
analytically in closed form
\begin{equation}
\label{eq8}
\begin{array}{l}
 c_{0k}  = \left. {\frac{{\sqrt \pi  {\mathop{\rm erf}\nolimits} \left( r \right) - 2are^{ - a^2 r^2 } }}{{4a^3 }}} \right|_{k}^{(k + 1)}  \\
 d_{jk}  = 2^{j/2} \left. {\frac{{\sqrt \pi  {\mathop{\rm erf}\nolimits} \left( r \right) - 2are^{ - a^2 r^2 } }}{{4a^3 }}} \right|_{2^{ - j} k}^{2^{ - j} (k + {\textstyle{1 \over 2}})}  -  \\
\\
 \quad \quad 2^{j/2} \left. {\frac{{\sqrt \pi  {\mathop{\rm erf}\nolimits} \left( r \right) - 2are^{ - a^2 r^2 } }}{{4a^3 }}} \right|_{2^{ - j} (k + {\textstyle{1 \over 2}})}^{2^{ - j} (k + 1)}  \\
 \end{array}
\end{equation}

\begin{figure}
\label{fig1}
\includegraphics[width=\textwidth]{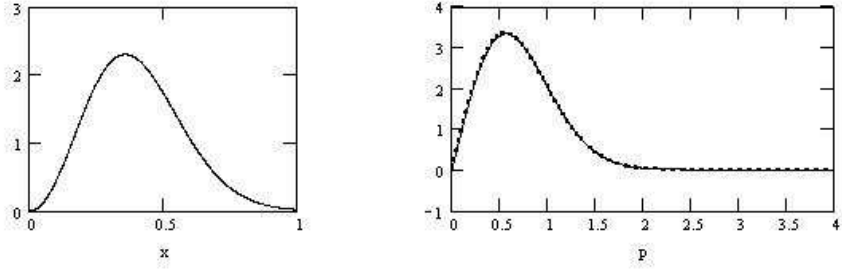}
  \caption{Original function (left) and Hankel fransform (right)}
\end{figure}

\begin{figure}
\label{f2}
\includegraphics[width=\textwidth]{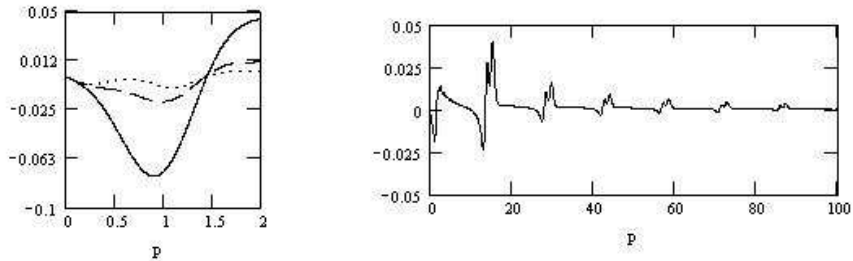}
  \caption{Transform's error}
\end{figure}

    Thus (\ref{eq6}) with the coefficients (\ref{eq8}) is the exact representation of the
Hankel transform. Let us consider the approximate solution. Suppose
 the function is known (\ref{eq7}) only in the segment $[0,h]$.
Then there is  the series instead of (\ref{eq4}):
   \begin{equation}
\label{eq9}
 g(r) = {c_{0} \varphi _0 (r)
 + \sum\limits_{j = 0}^J  {\sum\limits_{k = 0}^{2^J-1}  {d_{jk}
\psi _{jk}^{} (r)} } }.
\end{equation}

If $J \to \infty$ then (\ref{eq9}) is exact for this truncated function.
 But practically
 we only use several first levels.
  For example we can see  the original function
(the replacement $r$ to $x=r/h$ is used) (\ref{eq7}) and the
transform at the Fig.1. One can see that exact transform (solid
line)
 and the transform at level $J=3$ (dotted line)
 coincide at this figure.
The absolute errors between the exact transform and the
approximate transform at the levels $J=2$ (solid line), $J=3$
(dashed line), and $J=4$ (dotted line) are represented at the
Fig.2 (left).
  It's oblivious that the error
is small in comparison with the values of the $F_1(p)$.  The
absolute error at the level $J=3$ in a wide range of $p$ is
plotted in the right side of the Fig.2. One can see that this
error has quasi-periodic oscillations because the function is
truncated. But they decrease with the growth of $p$ (and $J$) when
oscillations in classical Fast Hankel Transform~\cite{Astro}
increase.


\begin{thebibliography}{99}
\bibitem{MMF}
 J. Mathews, R.L. Walker Mathematical methods of physics.
 New-York- Amsterdam, W.A. Benjamin, Inc, 1964
\bibitem{SF}
M. Abramowitz and I. A. Stegun (eds.), Handbook of mathematical functions with
formulas, graphs and mathematical tables, National Bureau of Standards Applied
Mathematics Series, vol. 55, U. S. Government Printing Office, Washington, D. C.,
1964.
\bibitem{Geo}
J. Zhao, W.W.M. Dai, S. Kapur, D.E. Long Efficient Three-Dimensional Extraction
Based on Statuic and Full-Wave Layered Green's Function. Proceedings of the 35
Conference on Design Automation, Moscon center, San Francisco, USA, June
15-19, 1998. ACM Press, 1998, pp. 224-229
\bibitem{Astro}
Hamilton A.J.S. Uncorrelated modes of the Nonlinear Power
Spectrum. Mon. Not. Roy. Astron. Soc. 2000, Iss. 312, pp. 257-284.
\bibitem{LO}
D. W. Lozier and F. W. J. Olver, Numerical evaluation of special functions, Mathematics
of Computation 1943-1993: A Half-Century of Computational Mathematics
(W. Gautschi, ed.), Proceedings of Symposia in Applied Mathematics, vol. 48, American
Mathematical Society, Providence, Rhode Island 02940, 1994, pp. 79-125. See
also http://math.nist.gov/nesf/.
\bibitem{FHT}
A. E. Siegman, Quasi fast Hankel transform, Optics Lett. 1 (1977), pp. 13-15.
\bibitem{Filon}
R. Barakat and E. Parshall, Numerical evaluation of the zero-order Hankel transform
using Filon quadrature philosophy, Appl. Math. Lett. 9 (1996), no. 5, pp. 21-26.
\bibitem{Chui}
C.K. Chui An introduction to Wavelets. Academic Press, 1992.
\end{thebibliography}
\end{document}